\documentclass{article}
\usepackage[utf8]{inputenc}
\usepackage{graphicx}
 
\usepackage{graphicx}
\usepackage{subfigure}
\usepackage{amsmath}
\usepackage{xcolor}
\usepackage{booktabs}
\usepackage{color, soul}
\usepackage{amssymb,latexsym}
\usepackage{amsmath,amsthm}
\usepackage{enumerate}
\usepackage{caption}
\usepackage[labelfont=bf]{caption}
\usepackage[labelsep=period]{caption}
\usepackage[mathscr]{eucal}
\usepackage{float}
\floatstyle{plaintop}
\restylefloat{table}
\usepackage[nobysame]{amsrefs}

\restylefloat{table}
\usepackage{placeins}
\usepackage{graphics} 
\usepackage{graphicx}
\graphicspath{ {./images/} }
\usepackage{float}
\title{Nonlocal theory of elasticity in the natural vibrations of curved CNTs with defects}
\author{Jaan Lellep and Shahid Mubasshar}
\date{}

\begin{document}

\maketitle

\section*{Abstract}
The current study presents the natural vibrations of single-walled carbon nanotubes (SWCNTs) using the nonlocal theory of elasticity. The study will help to develop nanodevices. The study focuses on the natural vibrations of curved or arch-like carbon nanotubes (CNTs). It is considered that the CNT has a crack-like defect and is simply supported (SS) at both ends. To model the equations of motion for the CNTs, a curved nanobeam approach is adopted, and the governing equations are solved using the separation of variables method. This research investigates the natural frequencies of three types of CNTs: armchair, zigzag, and chiral. A qualitative validation study demonstrates that the obtained results align with those published in the literature. Notably, the study reveals that small-scale effects and defects in the nanotubes influence the natural frequencies of CNTs.
\section{Introduction}
Nanostructures, such as nanotubes, nanoplates, nanoshells, nanoarches, and nano-cones, have exceptional mechanical and thermal properties due to their nanoscale size \cite{dai1996nanotubes}. The discovery of CNTs by Ijima in 1991 led to the revelation of numerous applications of these nanotubes. CNTs serve as fundamental components in the construction of nanodevices, nanosensors and nanoscillators \cites{gibson2007vibrations, journet2012carbon}. Consequently, researchers have extensively studied the dynamic and static behaviour of CNTs experimentally and theoretically \cites{lim2006coherent, lim2009coherent, babic2003intrinsic}.
Initially, vibration in CNTs was studied using classical continuum mechanics \cites{demir2010free, wang2006timoshenko}, neglecting small scale effects such as electric forces and Van der Waals forces. However, experimental studies have shown that these small scale effects, including the lattice spacing between individual atoms, may not be neglected since the material will not be homogenized into a continuum \cite{behera2017recent}. Consequently, various theories have been developed to investigate the dynamic behaviour of nanomaterials. One widely accepted approach is the nonlocal theory of elasticity, formulated by Eringen and his coworkers \cites{eringen1972nonlocal, eringen1983differential}, which has been extensively applied to study the vibrations of nanotubes. Subsequently, the nonlocal theory has also been employed to examine bending and buckling problems in nanomaterials \cites{reddy2008nonlocal, wang2006small}. Moreover, researchers have studied the effects of cracks on the static and dynamic behaviour of nanostructures \cites{hossain2020effect, mubasshar2022natural}, as cracks can form on the surfaces of nanomaterials due to environmental factors, improper handling and wear and tear. This study focuses on the natural frequencies of armchair, zigzag and chiral shape CNTs. The effects of cracks, nonlocal and other physical parameters on the natural frequency of CNTs are discussed.

\section{Formulation and solution of the problem}
\begin{figure}[H]
\centerline{\includegraphics[width=10cm,height=7cm]{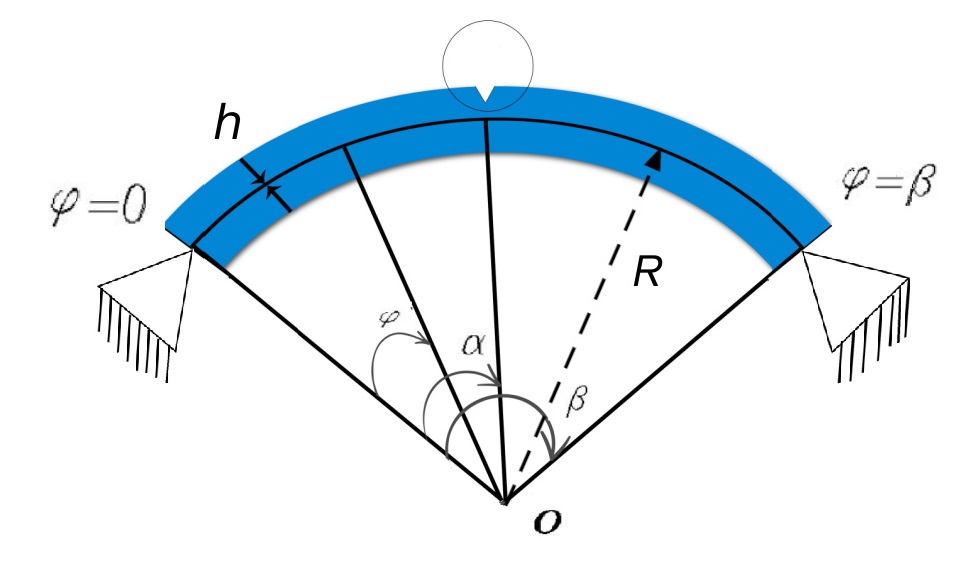}}
\caption{A schematic diagram of the simply supported curved CNT.}
\centering
\label{Fig.1}
\end{figure}

Creating an armchair, zigzag, and chiral CNTs from a graphene sheet involves rolling and manipulating the sheet differently, as shown in \cite{hussain2019nonlocal}. 
The geometry of a SS carbon nanotube (CNT) with a defect is shown in Fig.\ref{Fig.1}. Here $h$ is the thickness and $R$ is the radius of the CNT. The current cross-section is defined by $\varphi$, and the edges of the arch correspond to $\varphi=0$ and $\varphi=\beta$. It is assumed that a stable crack whose length is $c<h$ is located at $\varphi=\alpha$. It is also assumed that the cross-section of the CNT is circular. 
 The equilibrium equations for the curved beam element are established by \cite{soedel2004vibrations} and have been extensively validated through experimental and numerical studies in the analysis of vibration of nanostructures \cite{lellep2022free}. According to this theory, the equilibrium equations are written as 
\begin{equation}
\label{eqn:1}
\begin{aligned}
  \frac {\partial{N}} {\partial s}{}+{} \frac {Q}{R}&+p_{u}&= &&\bar \rho h\frac{\partial^2 U}{\partial t^2},\\
 \frac {\partial{Q}} {\partial s} {}-{}\frac {N}{R}&+p&= &&\bar \rho h\frac{\partial^2 W}{\partial t^2},\\
\frac {\partial{M}} {\partial{ s}}&-Q&=&&0,
\end{aligned}
\end{equation}

where $N$, $M$, and $Q$ represent the membrane force, bending moment, and shear force, respectively. The length $s$ of the nanoarch is related to the angle $\varphi$ by $s = R\varphi$.
In (\ref{eqn:1}), $U$ and $W$ correspond to the displacements in the circumferential and transverse directions, respectively, while $p_{u}$ and $p$ denote external loads in those directions. The variable $h$ represents the CNT's thickness, $\bar{\rho}$ is the mass per unit length, and $b$ is the width. Strain-displacement equations, as stated in \cite{soedel2004vibrations} are
\begin{equation}
\label{eqn:3}
\varepsilon=\frac {1} {R}W+\frac {\partial{U}} {\partial {s}}
\end{equation}
and
\begin{equation}
\label{eqn:4}
\varkappa=-\frac {\partial^2{U}} {\partial s^2}+\frac {\partial{U}} {\partial{s}}\frac{1}{R}.
\end{equation}
Here, $\varepsilon$ denotes the relative extension of the curved element and $\varkappa$ represents the curvature of the nanotube's middle line.
The classical theory of elasticity describes Hooke's law as
\begin{equation}
M_{c} = EI\varkappa, 
\end{equation}
where $M_{c}$ is the bending moment, $E$ stands for Young's modulus, and $I=\frac{\pi d^{4}}{64}$ denotes the moment of inertia with diameter $d$.  The material of the CNT adheres to the constitutive equations of the nonlocal theory of elasticity \cites{behera2017recent, wang2006small}.  According to the theory, the constitutive equation is given as
\begin{equation}
\label{eqn:7}
\sigma_{ij} - \eta\nabla^{2}\sigma_{ij} = \sigma^c_{ij}, 
\end{equation}
where $\sigma^c_{ij}$ is nonlocal stress tensor, $\nabla$ represents the Laplacian operator, and $\eta$ is a material constant defined as $\eta = (e_{0}a)^2$, with $e_{0}$ being a calibration constant for different materials and $a$ is a characteristic internal length.
In this problem, we assume that $p_{u} = 0$, $U(\varphi, t)$ is small, and $\ddot{U}$ is also negligible. Under these conditions, (\ref{eqn:1}) can be simplified as
\begin{equation}
N^{'} = -Q, 
\end{equation}
\begin{equation}
Q^{'} = N + R(\rho h \ddot{W} - p), 
\end{equation}
\begin{equation}
M^{'} = RQ. 
\end{equation}

Here, prims denote differentiation with respect to the current angle $\varphi$. Furthermore, if $M(\varphi)$ and $N(\varphi)$ vanish at the same cross-section, one can write
\begin{equation}
\label{eqn:9}
M = -RN. 
\end{equation}
To solve the governing equations, the method of separation of variables is applied, which assumes that 
 \begin{equation}
\label{eqn:16}
  W(\varphi, t)=X(\varphi) \sin(\omega t),
\end{equation}
where $X(\varphi)$ is a function of the coordinate $\varphi$.
Equations (\ref{eqn:9}, \ref{eqn:16}) can be written as 
\begin{equation}
\label{eqn:19}
X^{''''}+(2+K\eta)X^{''}+(1-K)X=0,
\end{equation}

where
\begin{equation}\nonumber
K=\frac{\omega ^{2} \rho h R^{4}}{EI},
\end{equation}
Cracks and other defects cause additional structural compliance, which can be calculated by using the method developed in \cite{dimarogonas1996vibration} and \cite{chondros2001vibration} and has been applied by many authors in the study of vibration problems \cites{mubasshar2022natural, hussain2019nonlocal, lellep2022free}. 
\section{Results and discussion}
\begin{figure}[H]
\centerline{\includegraphics[width=12cm,height=7cm]{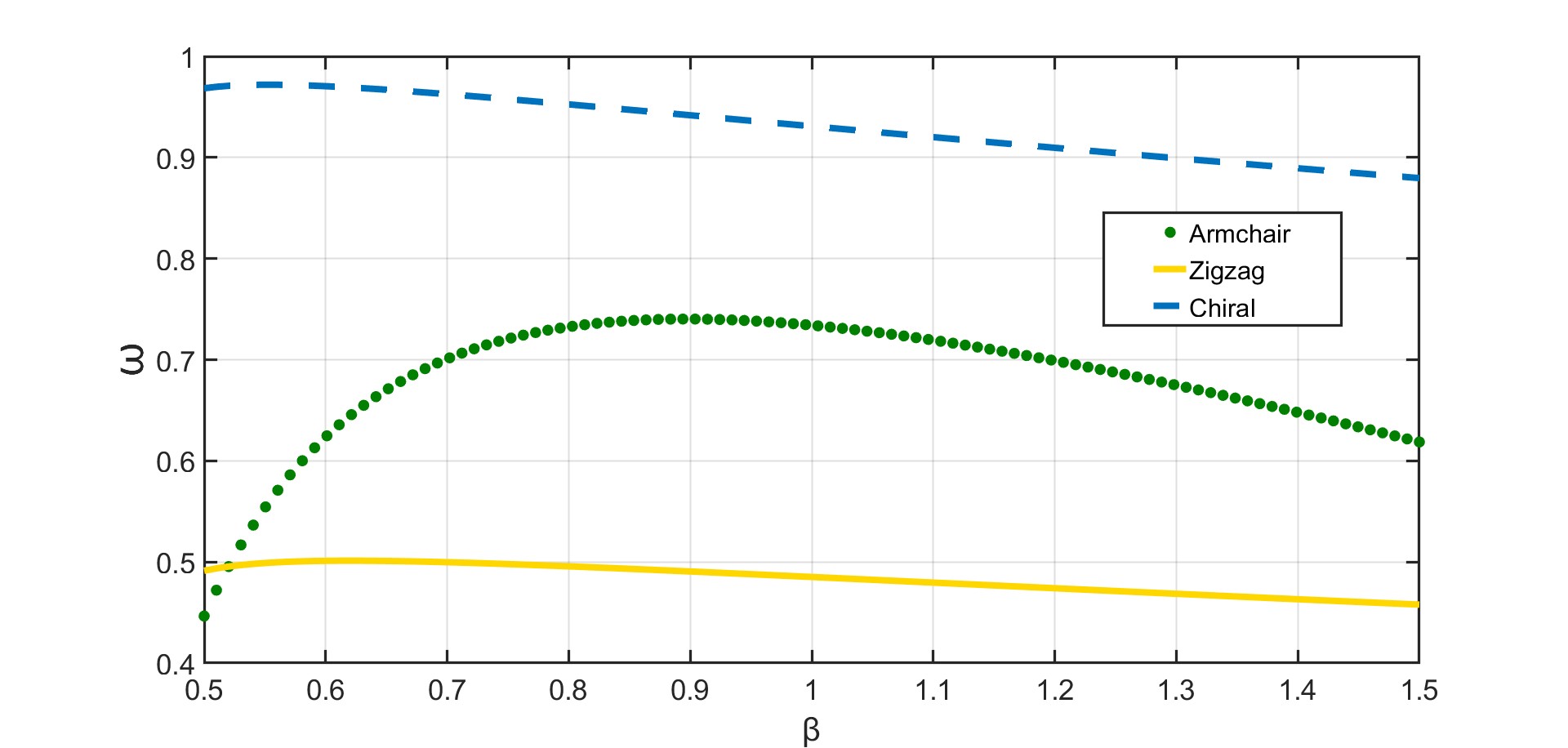}}
\vspace*{8pt}
\caption{Natural frequency versus central angle of the CNTs.}
\label{fig:fig.3}
\centering
\end{figure}

\begin{figure}[H]
\centerline{\includegraphics[width=12cm,height=7cm]{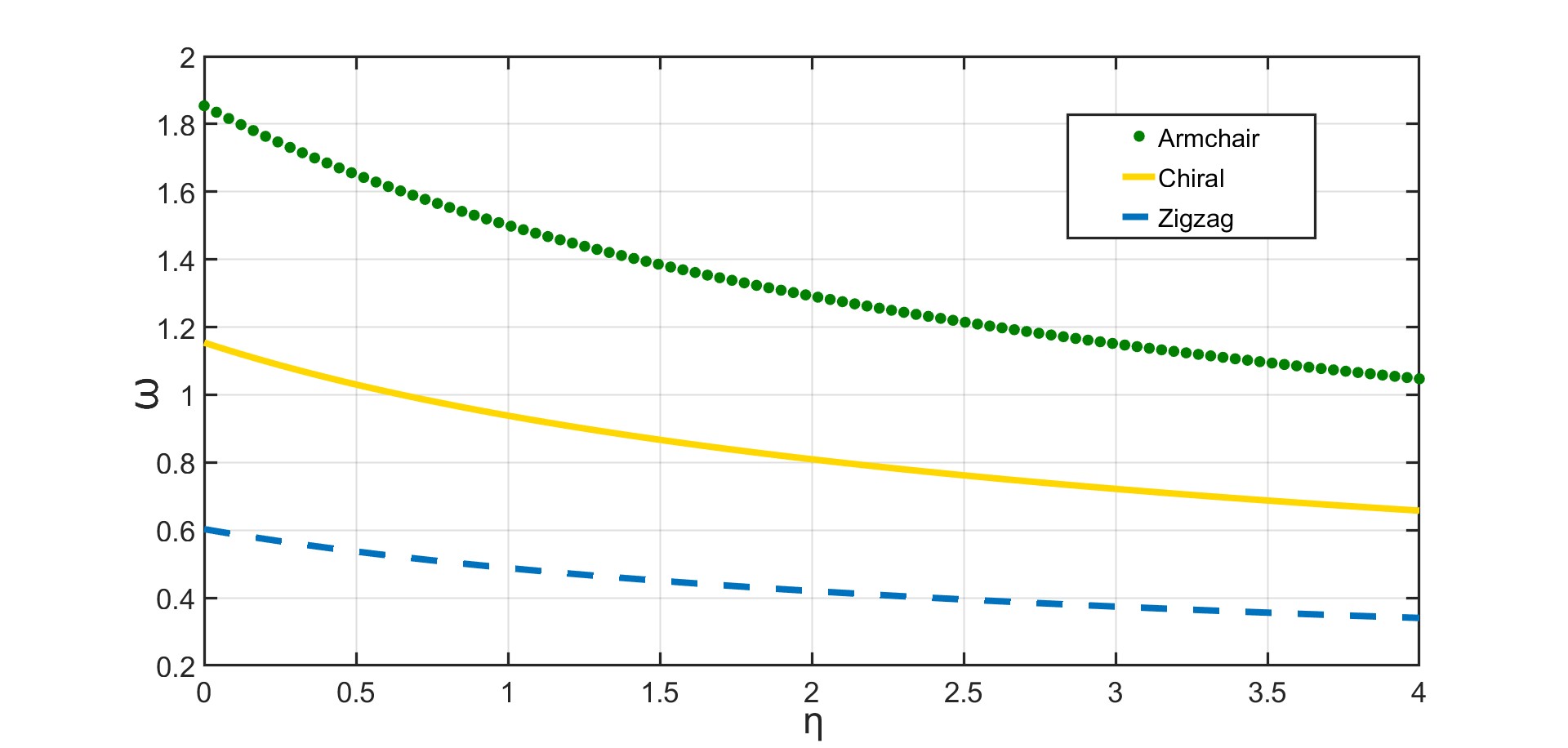}}
\vspace*{8pt}
\caption{Natural frequency versus the nonlocal parameter.}
\label{fig:fig.4}
\centering
\end{figure}

\begin{figure}[H]
\centerline{\includegraphics[width=12cm,height=7cm]{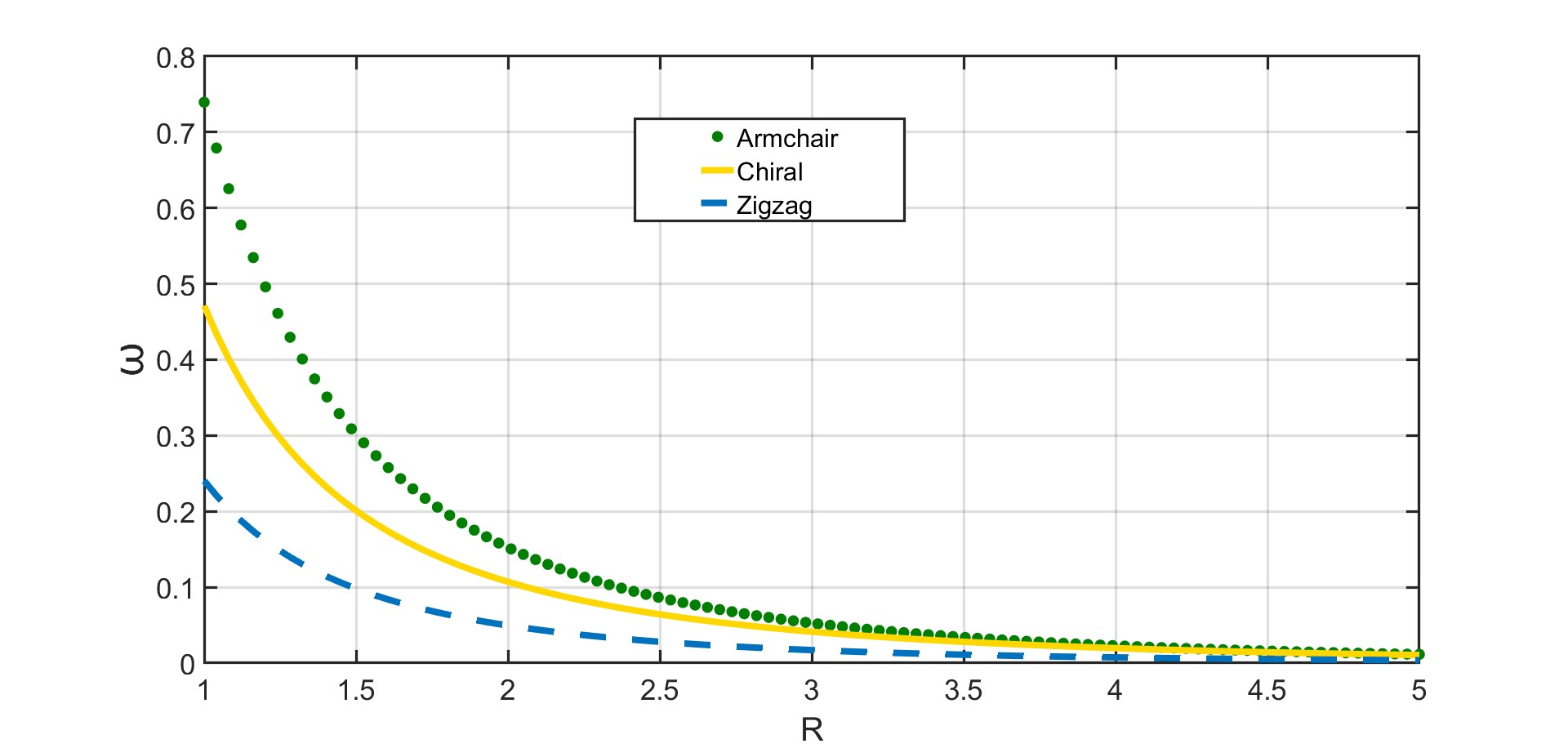}}
\vspace*{8pt}
\caption{Natural frequency versus radius of the CNTs.}
\label{fig:fig.5}
\centering
\end{figure}

\begin{figure}[H]
\centerline{\includegraphics[width=12cm,height=7cm]{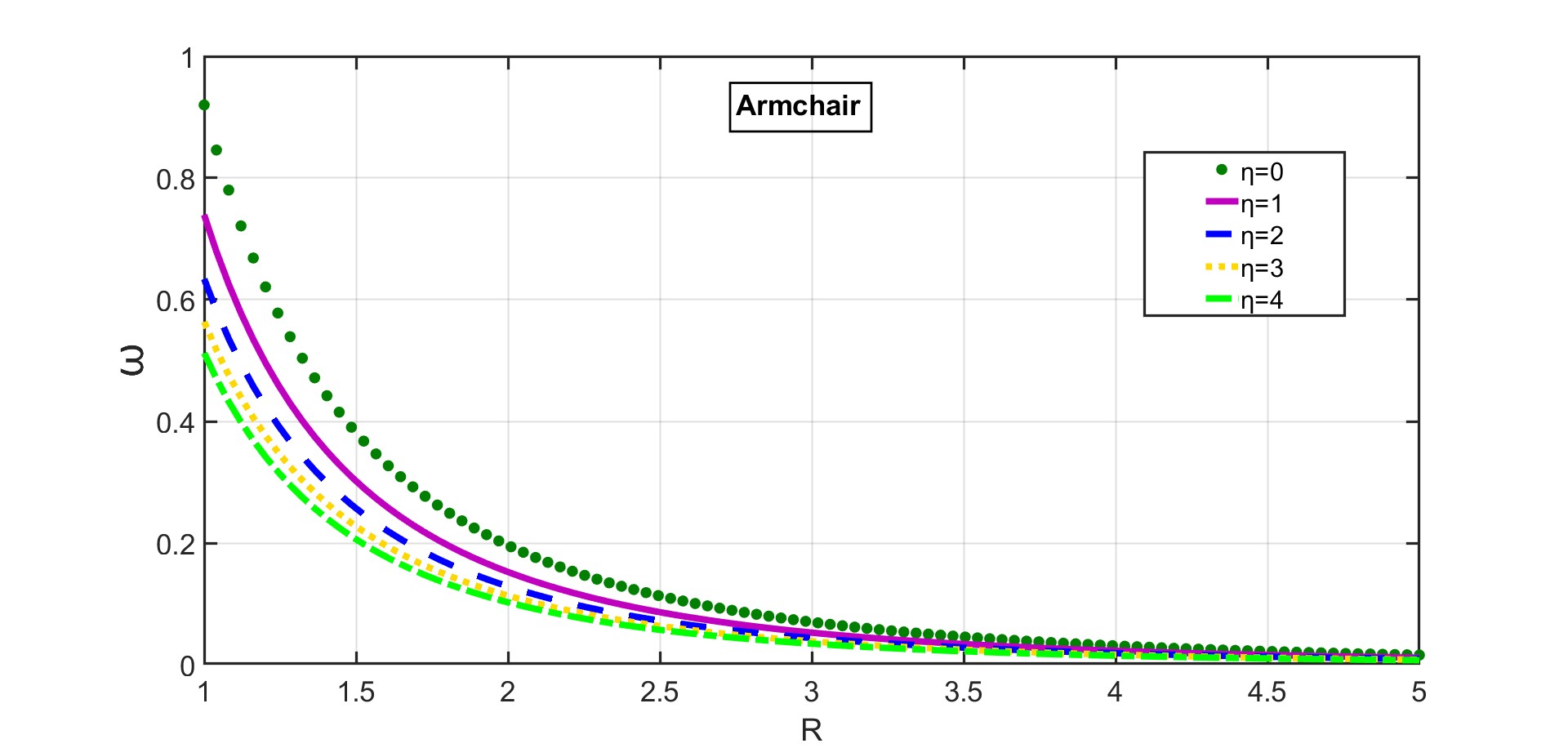}}
\vspace*{8pt}
\caption{Natural frequency versus radius of the armchair CNT.}
\label{fig:fig.6}
\centering
\end{figure}

\begin{figure}[H]
\centerline{\includegraphics[width=12cm,height=7cm]{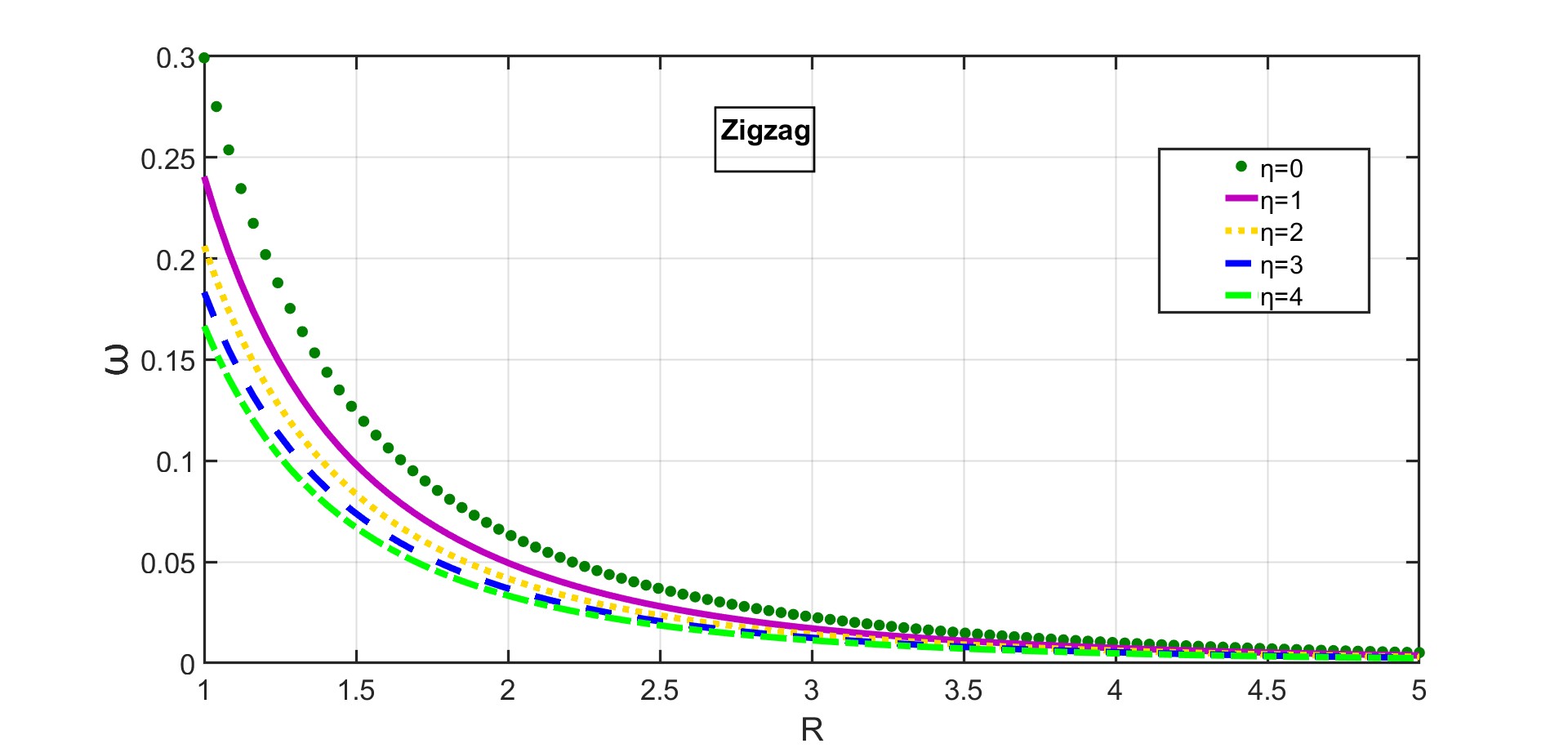}}
\vspace*{8pt}
\caption{Natural frequency versus radius of the zigzag CNT.}
\label{fig:fig.7}
\centering
\end{figure}

Numerical results are presented for the SWCNTs of different types with simply supported boundary conditions. Physical properties of the CNTs are taken from published literature \cite{gupta2010wall} and \cite{majeed2019vibration}. $\beta=1$, $\eta=1$ and crack length $c=4$ are considered where values of these parameters are not specified in the legends.  
In Fig. \ref{fig:fig.3}, the natural frequency of the CNTs is plotted against the central angle $\beta$. From the figure, one can observe that the natural frequency initially increases and then starts decreasing as we increase the $\beta$ value. It is also noticeable that the frequency for the chiral nanotubes is higher than that of armchair and zigzag nanotubes.
Fig. \ref{fig:fig.4} shows the natural frequency plotted against the nonlocal parameter $\eta$. The figure illustrates that the frequency consistently decreases as we increase the $\eta$. Furthermore, it can be seen that the zigzag nanotubes exhibit the lowest frequency compared to other nanotubes.
In Fig. \ref{fig:fig.5}, the natural frequency is plotted against the radius of the CNTs. It is clear from the figure that the frequency decreases as the value of the radius increases.
Figs. \ref{fig:fig.5}-\ref{fig:fig.7} depict the natural frequencies of armchair, zigzag, and chiral CNTs, plotted against the radius of nanotubes. Different curves in the figures correspond to different values of $\eta$. In Table \ref{tab:thai1}, the results are compared with \cite{thai2012nonlocal} when there is no crack in the nanotubes and central angle $\beta$ is very small. The table shows that the results are close to each other for different values of nonlocal parameter $\eta$.

\begin{figure}[H]
\centerline{\includegraphics[width=12cm,height=7cm]{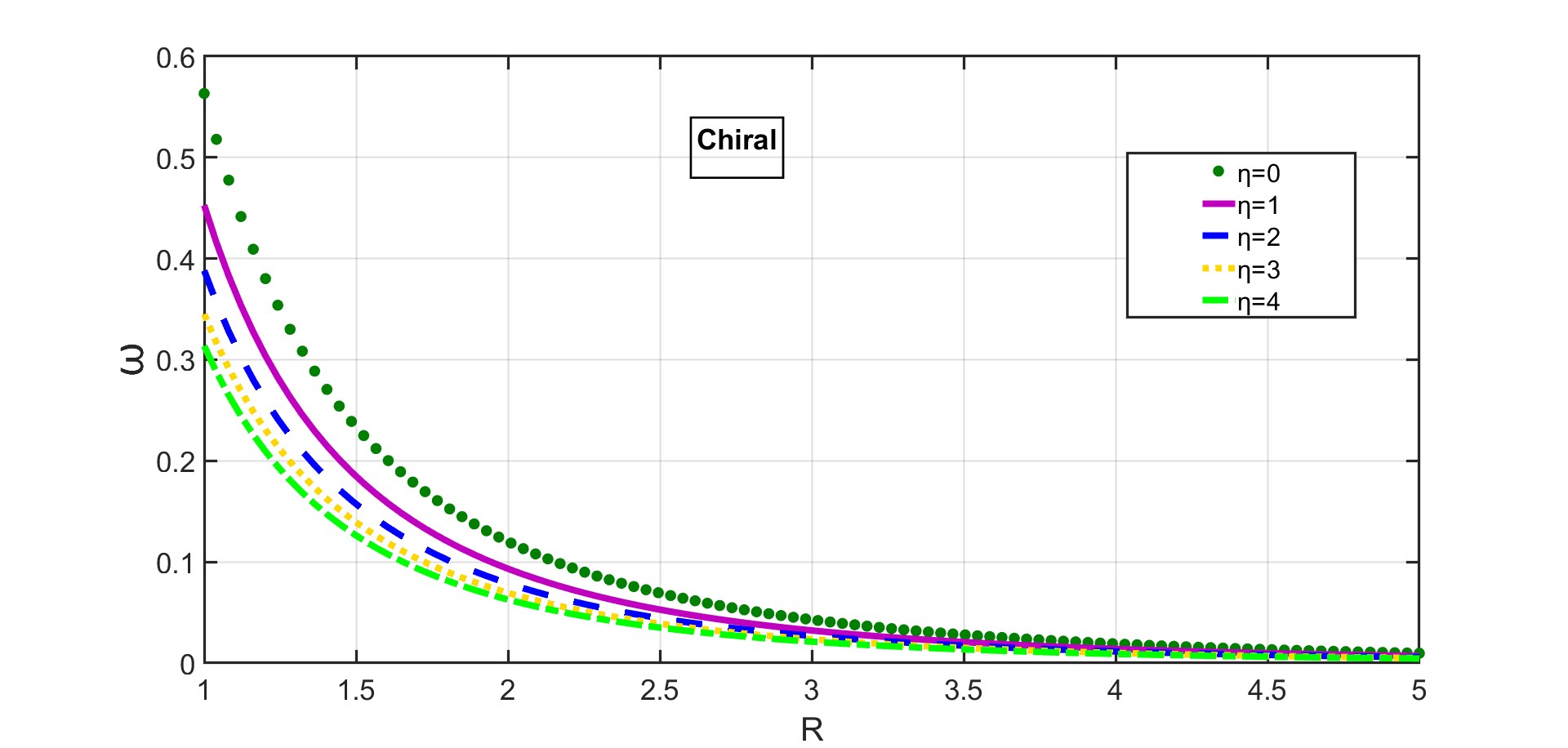}}
\vspace*{8pt}
\caption{Natural frequency versus radius of the chiral CNT.}
\label{fig:F8}
\centering
\end{figure}
\FloatBarrier
\begin{table}[hbt!]
\caption{Natural vibration of the nanobeams of the constant thickness.}
\centering
\begin{tabular}{ |c|c|c|c|c| } 
 \hline
 Mode&$\eta$&  Present & Thai \cite{thai2012nonlocal} \\ 
  \hline
 1&0& 9.75821 & 9.2745 \\ 
  
& 1 & 7.05584 & 8.8482 \\ 
 
 &2& 5.80188 & 8.4757 \\ 
 
 &3&  5.04192 & 8.1466 \\ 
 
 &4&  4.51883 & 7.8530 \\ 
 \hline
\end{tabular}

\label{tab:thai1}
\end{table}
\FloatBarrier

\section{Conclusion}
Natural frequencies of CNTs of three different shapes, armchair, zigzag and chiral, are calculated based on the nonlocal theory of elasticity. The small scale effects on the natural frequency of CNTs are exhibited. The effects of the other physical parameters, such as radius and central angle on the natural frequency are also calculated. The results show that the natural frequency of the CNTs decreases as we increase the nonlocality and the radius of the curved shape CNTs. To verify the methodology, the obtained results are compared with those available in the literature. The study is helpful for engineers in the field of the construction of nanodevices.   

\bibliographystyle{plain}
\bibliography{main}
\end{document}